\documentclass[12pt]{article}

\usepackage{amsmath}
\usepackage{amssymb}
\usepackage{microtype}
\usepackage{amsthm}
\usepackage[margin=33mm]{geometry}

\newtheorem{thm}{Theorem}[section]
\newtheorem{defn}[thm]{Definition}

\newtheorem{conj}[thm]{Conjecture}

\begin{document}

\renewcommand{\thefootnote}{\fnsymbol{footnote}}

\title{Stirling numbers of forests and cycles}

\author{David Galvin\thanks{dgalvin1@nd.edu; Department of Mathematics,
University of Notre Dame, Notre Dame IN 46556. Supported in part by National Security Agency grant H98230-10-1-0364.}
\and Do Trong Thanh\thanks{dtrongth@nd.edu; University of Notre Dame, Notre Dame IN 46556. Supported by Glynn Family Honors Program.}}

\date{\today}

\maketitle

\begin{abstract}
For a graph $G$ and a positive integer $k$, the {\em graphical Stirling number} $S(G,k)$ is the number of partitions of the vertex set of $G$ into $k$ non-empty independent sets. Equivalently it is the number of proper colorings of $G$ that use exactly $k$ colors, with two colorings identified if they differ only on the names of the colors. If $G$ is the empty graph on $n$ vertices then $S(G,k)$ reduces to  $S(n,k)$, the familiar Stirling number  of the second kind.

In this note we first consider Stirling numbers of forests. We show that if $(F^{c(n)}_n)_{n\geq 0}$ is any sequence of forests with $F^{c(n)}_n$ having $n$ vertices and $c(n)=o(\sqrt{n/\log n})$ components, and if $X^{c(n)}_n$ is a random variable that takes value $k$ with probability proportional to $S(F^{c(n)}_n,k)$ (that is, $X^{c(n)}_n$ is the number of classes in a uniformly chosen partition of $F^{c(n)}_n$ into non-empty independent sets), then $X^{c(n)}_n$ is asymptotically normal, meaning that suitably normalized it tends in distribution to the standard normal. This generalizes a seminal result of Harper on the ordinary Stirling numbers. Along the way we give recurrences for calculating the generating functions of the sequences $(S(F^c_n,k))_{k \geq 0}$, show that these functions have all real zeroes, and exhibit three different interlacing patterns between the zeroes of pairs of consecutive generating functions.

We next consider Stirling numbers of cycles. We establish asymptotic normality for the number of classes in a uniformly chosen partition of $C_n$ (the cycle on $n$ vertices) into non-empty independent sets. We give a recurrence for calculating the generating function of the sequence $(S(C_n,k))_{k \geq 0}$, and use this to give a direct proof of a log-concavity result that had previously only been arrived at in a very indirect way.
\end{abstract}

\section{Introduction}

Let $G=(V,E)$ be a (simple, finite, loopless) graph on $n$ vertices. An {\em independent set} in $G$ is a subset of the vertices, no two of which are adjacent. For each integer $k$ we set
$$
S(G,k) = \left|\{\mbox{partitions of $V$ into $k$ non-empty independent sets}\}\right|.
$$
We may interpret $S(G,k)$ as the number of proper $k$-colorings of $G$ that use all $k$ colors, with two colorings identified if they are identical up to the names of the colors. As far as we can discover, this parameter was first explicitly considered by Tomescu \cite{Tomescu}. When $G=E_n$, the graph with no edges, $S(G,k)$ is just the Stirling number of the second kind $S(n,k)$, the number of partitions of a set of size $n$ into exactly $k$ non-empty classes. It is for this reason that in \cite{DuncanPeele}, Duncan and Peele refer to $S(G,k)$ as a {\em graphical Stirling number}. Motivated by the connection to colorings, Goldman, Joichi and White \cite{GoldmanJoichiWhite} refer to the sequence $(S(G,k))_{k\geq 0}$ as the {\em chromatic vector} of $G$.

Recall that associated with each $G$ there is a polynomial $\chi_G(x)$ of degree $n$, the {\em chromatic polynomial} of $G$, whose value at each positive integer $x$ is the number of proper $x$-colorings of $G$, that is, the number of functions $f:V\rightarrow \{1,\ldots, x\}$ satisfying $f(u)\neq f(v)$ whenever $uv \in E$. The chromatic polynomial of $G$ determines the sequence $(S(G,k))_{k \geq 0}$, and vice-versa. On the one hand inclusion-exclusion gives
\begin{equation} \label{eq-chrom}
S(G,k) = \frac{1}{k!}\sum_{i=0}^k (-1)^i{k \choose i}\chi_G(k-i),
\end{equation}
while on the other hand
$$
\chi_G(x) = \sum_{k \geq 0} S(G,k)x_{(k)}
$$
where $x_{(k)} = x(x-1) \ldots (x-k+1)$. Indeed, given a palette of $x$ colors, for each $k$ there are $S(G,k)$ ways to partition the vertex set into $k$ non-empty color classes, and $x_{(k)}$ ways to assign colors the classes.

Note that $S(G,k)=0$ whenever $k < \chi(G)$ (where $\chi(G)$ is the chromatic number of $G$, the least positive integer $x$ for which $\chi_G(x)>0$) and also whenever $k > n$. Also, $S(G,k)$ is non-zero at $k=\chi(G)$ and takes value $1$ at $k=n$. It follows that the generating function of the sequence $(S(G,k))_{k \geq 0}$,
$$
\sigma(G,x) = \sum_{k \geq 0} S(G,k)x^k
$$
is a monic polynomial of degree $n$, and has a zero with multiplicity $\chi(G)$ at $0$. We refer to $\sigma(G,x)$ as the {\em Stirling polynomial} of $G$.

The Stirling polynomial (arrived at via an alternate definition) was first considered by Korfhage \cite{Korfhage}, who established some basic properties and calculated it explicitly for some particular families of graphs (including forests). It was then extensively studied by Brenti \cite{Brenti2} and Brenti, Royle and Wagner \cite{BrentiRoyleWagner}, who established large classes of graphs (including cycles) for which the Stirling polynomial has all real zeroes. Graphical Stirling numbers (in various guises) for some particular families of graphs have also been considered by Farrell and Whitehead \cite{FarrellWhitehead}, Duncan and Peele \cite{DuncanPeele}, Mohr and Porter \cite{MohrPorter}, Goldman, Joichi and White \cite{GoldmanJoichiWhite}, and Munagi \cite{Munagi}.

There are some important consequences of the generating function of a sequence $(a_k)_{k=0}^n$ of non-negative terms having all real zeroes. A theorem of Newton (see for example \cite[page 504]{Stanley}) implies that such a sequence satisfies the ultra log-concavity relation
$$
\left(\frac{a_k}{{n \choose k}}\right)^2 \geq \left(\frac{a_{k-1}}{{n \choose k-1}}\right)\left(\frac{a_{k+1}}{{n \choose k+1}}\right)
$$
for $k=1, \ldots, n-1$, which in turn implies the usual log-concavity relation $a_k^2 \geq a_{k-1}a_{k+1}$. If the sequence has no internal zeros (i.e., no $k$ with $a_{k-1}\neq 0$, $a_k=0$ and $a_{k+1}\neq 0$) then log-concavity in turn implies that the sequence is unimodal, that is, that there is $0 \leq k \leq n$ with $a_0 \leq \ldots \leq a_k \geq \ldots \geq a_n$.

Another consequence of real zeroes is that if $X$ is a random variable supported on $\{0, \ldots, n\}$ that takes value $k$ with probability proportional to $a_k$, then $X$ may be realized as the sum of $n$ independent Bernoulli random variables, leading to the possibility of a central limit theorem and asymptotic normality.
\begin{defn}
A sequence $(X_n)_{n \geq c}$ is said to be {\em asymptotically normal} if, for each $x \in {\mathbb R}$, we have
$$
\lim_{n \rightarrow \infty} \Pr \left(\frac{X_n-E(X_n)}{\sqrt{{\rm Var}(X_n)}} \leq x\right) = \frac{1}{\sqrt{2\pi}} \int_{-\infty}^x e^{-x^2/2}~\!dx,
$$
with the convergence uniform in $x$.
\end{defn}
In other words, a sequence is asymptotically normal if, suitably normalized, it tends in distribution to the standard normal. Note that if $a_k$ counts the number of members of size $k$ of some family, then $X$ measures the size of a uniformly selected member of the family.

A seminal result in this area is due to Harper \cite{Harper}. He proved that for each $n$, the generating function of the sequence $(S(n,k))_{k \geq 0}$ of the Stirling numbers of the second kind has all real zeroes, and used this to prove that the sequence of random variables $(X^{\rm Stirling}_n)_{n \geq 0}$, where $X^{\rm Stirling}_n$ is the number of classes in a uniformly chosen partition of a set of size $n$ into non-empty classes, is asymptotically normal.
The first aim of this note is to obtain an analog of Harper's result (and, as it turns out, a generalization) for the Stirling numbers of forests. In what follows $B_n=\sum_{k=0}^n S(n,k)$ is the $n$th Bell number, and $W(n) = \log n -\log \log n + O(1)$ is the Lambert $W$ function, the unique positive solution to $W(n)e^{W(n)} = n$.
\begin{thm} \label{thm-forest-asy-nor}
Let $(F^{c(n)}_n)_{n \geq 0}$ be a sequence of forests, with $F^{c(n)}_n$ having $n$ vertices and $c(n)$ components. Let $X^{c(n)}_n$ be the number of classes in a uniformly chosen partition of the vertex set of $F^{c(n)}_n$ into non-empty independent sets. There is a constant $C>0$ such that if $c(n)<C\sqrt{n/\log n}$ for all sufficiently large $n$, then the sequence $(X^{c(n)}_n)_{n \geq 0}$ is asymptotically normal, with
\begin{eqnarray}
E(X^{c(n)}_n) & = & \frac{\sum_{i \geq 0} {c(n)-1 \choose i}B_{n-i}}{\sum_{i \geq 0} {c(n)-1 \choose i}B_{n-1-i}} \label{forest-exp} \\
& = & \frac{n}{W(n)} + O\left(\frac{1}{\log n}\right), \label{forest-exp-est}
\end{eqnarray}
and
\begin{eqnarray}
{\rm Var}(X^{c(n)}_n) & = & \frac{\sum_{i \geq 0} {c(n)-1 \choose i}B_{n+1-i}}{\sum_{i \geq 0} {c(n)-1 \choose i}B_{n-1-i}} -\left(\frac{\sum_{i \geq 0} {c(n)-1 \choose i}B_{n-i}}{\sum_{i \geq 0} {c(n)-1 \choose i}B_{n-1-i}}\right)^2 - 1 \label{forest-var} \\
& = & \frac{n}{W(n)(W(n)+1)} + O\left(\frac{c(n)^2}{\log n}\right). \label{forest-var-est}
\end{eqnarray}
\end{thm}
\begin{conj} \label{conj-forests}
The sequence $(X^{c(n)}_n)_{n \geq 0}$ is asymptotically normal for all $1 \leq c(n) \leq n$.
\end{conj}

The denominator of (\ref{forest-exp}) turns out to be $\sum_{k \geq 0} S(F_n^{c(n)},k)$, the total number of partitions of $F_n^{c(n)}$ into non-empty independent sets. For general $G$, the quantity $\sum_{k \geq 0} S(G,k)$ is referred to in \cite{DuncanPeele} as the {\em Bell number} of $G$. For $G=E_n$, the Bell number is just the ordinary Bell number $B_n$.

We begin the proof of Theorem \ref{thm-forest-asy-nor} by deriving an explicit expression for $S(F^c_n,k)$ in terms of ordinary Stirling numbers. The same expression was obtain by Korfhage \cite{Korfhage} using induction; our direct derivation is based on the fact that the chromatic polynomial of a forest depends only on how many  vertices and components it has. This expression allows us to obtain a number of recurrences for $\sigma(F^c_n,x)$, which allow us to establish that $\sigma(F^c_n,x)$ has all real zeros for all $n$ and $c$, (a result already implicit in \cite{Brenti2}) and that moreover the zeroes of the $\sigma(F^c_n,x)$'s exhibit a number of nice interlacing patterns, none of which have been previously observed.

\begin{defn}
Given two reals polynomials $f$ and $g$, we say that {\em $f \prec g$} if all of the following conditions hold:
\begin{enumerate}
\item the zeroes of both $f$ and $g$ are all real and non-positive, and all negative zeroes of both have multiplicity one;
\item the number of negative zeroes of $g$ is either the same or one greater than the number of negative zeroes of $f$; and
\item if $x_1, \ldots, x_m$ are the negative zeroes of $g$ and $y_1, \ldots, y_n$ are the negative zeroes of $f$, both written in decreasing order, then
    $$
    x_1 > y_1 > x_2 > y_2 > \ldots.
    $$
\end{enumerate}
\end{defn}
Note that this is not the same as saying that the roots of $f$ and $g$ interlace in the normal sense; we introduce the notation $f \prec g$ to deal with the fact that most of the polynomials we will encounter have multiple zeroes at $0$.
\begin{thm} \label{thm-interlacing}
Fix $c \geq 1$ and $n \geq c$. We have the following interlacing relations between the zeroes of the $\sigma(F^c_n,x)$'s.
\begin{enumerate}
\item $\sigma(F^c_{c+1},x) \prec \sigma(F^c_c,x)$
\item $\sigma(F^c_n,x) \prec \sigma(F^c_{n+1},x)$ for all $n \geq c+1$
\item $\sigma(F^c_n,x) \prec \sigma(F^{c+1}_{n+1},x)$
\item $\sigma(F^c_{c+1},x) \prec \sigma(F^{c+1}_{c+1},x)$
\item $\sigma(F^{c+1}_{n+1},x) \prec \sigma(F^c_{n+1},x)$ for all $n \geq c+1$
\end{enumerate}
\end{thm}

$$
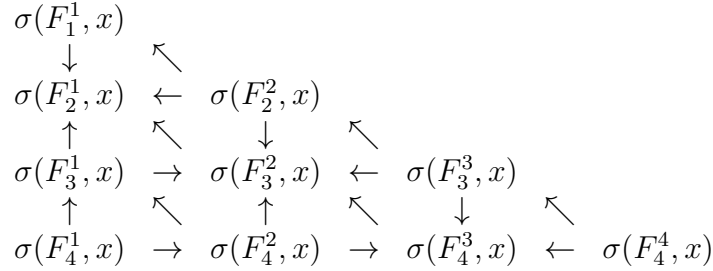
\begin{array}{c}
\begin{array}{ccccccc}
\sigma(F^1_1,x) &   &   &   &   &   &    \\
\downarrow      &  \nwarrow &   &   &   &   &  \\
\sigma(F^1_2,x) & \leftarrow & \sigma(F^2_2,x) &   &   &   &  \\
\uparrow                 &  \nwarrow &  \downarrow & \nwarrow  &   &   &  \\
\sigma(F^1_3,x) & \rightarrow &  \sigma(F^2_3,x) &  \leftarrow &  \sigma(F^3_3,x) &  & \\
\uparrow  & \nwarrow  & \uparrow  & \nwarrow  &  \downarrow & \nwarrow  &  \\
\sigma(F^1_4,x) & \rightarrow  & \sigma(F^2_4,x)  & \rightarrow  &  \sigma(F^3_4,x) & \leftarrow  & \sigma(F^4_4,x)
\end{array} \\
\\
\mbox{{\small The interlacing patterns between the $\sigma(F^c_n,x)'s$.}} \\
\mbox{{\small An arrow from $P$ to $Q$ (with $P$ at the head) indicates that $P \prec Q$.}}
\end{array}
$$

We prove Theorem \ref{thm-interlacing} in Section \ref{sec-forest}, where we also derive (\ref{forest-exp}) and (\ref{forest-var}). Known asymptotic estimates for Bell numbers, together with a general theory of asymptotic normality (which we briefly discuss in Section \ref{sec-asy-nor}),
then allows us to complete the proof of Theorem \ref{thm-forest-asy-nor}. The details are in Section \ref{sec-bell}.

The second aim of this note is to obtain an analog of Harper's result for the Stirling numbers of cycles.
\begin{thm} \label{thm-cycles-asy-nor}
Let $X^{\rm cycle}_n$ be the number of classes in a uniformly chosen partition of the vertex set of the cycle $C_n$ on $n$ vertices into non-empty independent sets. The sequence $(X^{\rm cycle}_n)_{n \geq 0}$ is asymptotically normal,
with
\begin{eqnarray}
E(X^{\rm cycle}_n) & = & \frac{\sum_{i \geq 0} (-1)^i B_{n-i}}{\sum_{i \geq 0} (-1)^i B_{n-1-i}} \label{cycle-exp} \\
&  = & \frac{n}{W(n)} + O\left(\frac{1}{\log n}\right) \label{cycle-exp-est}
\end{eqnarray}
and
\begin{eqnarray}
{\rm Var}(X_n) & = & \frac{\sum_{i \geq 0} (-1)^i B_{n+1-i}}{\sum_{i \geq 0} (-1)^i B_{n-1-i}} - \left(\frac{\sum_{i \geq 0} (-1)^i B_{n-i}}{\sum_{i \geq 0} (-1)^i B_{n-1-i}}\right)^2 - 1 \label{cycle-var} \\
& = & \frac{n}{W(n)(W(n)+1)} + O\left(\frac{1}{\log n}\right). \label{cycle-var-est}
\end{eqnarray}
\end{thm}

In the course of the proof of Theorem \ref{thm-forest-asy-nor} we will see that the Bell number of the cycle $C_n$ is $\sum_{i \geq 0} (-1)^i B_{n-1-i}$. This expression also occurs in \cite{CzabarkaErdosJohnsonKupczokSzekely}, where it shown to be the number of partitions of a set of size $n$ into non-empty classes each of size at least $2$. One of the main concerns of \cite{CzabarkaErdosJohnsonKupczokSzekely} is the quantity $S^\star(n,k)$, the number of such partitions with exactly $k$ classes, and it is shown that the random variable taking value $k$ with probability proportional to $S^\star(n,k)$ is asymptotically normal. We do not see any connection, though, between this question and the question of graphical Stirling numbers of cycles, other than the coincidence between the Bell number of $C_n$ and the total number of singleton-free partitions of $[n]$.

As with Theorem \ref{thm-forest-asy-nor}, we begin by deriving a recurrence for $\sigma(C_n,x)$. Ideally we would like to use this recurrence to give a direct proof of the fact that for each $n \geq 3$, $\sigma(C_n,x)$ has all real zeroes. This was already established in \cite{BrentiRoyleWagner} (after having been conjectured by Brenti in \cite{Brenti2}). The methods of \cite{BrentiRoyleWagner} are circuitous, involving results of Wagner on partition polynomials \cite{Wagner}, a detailed study of an operator that converts the Stirling polynomial into the chromatic polynomial, and facts about the locations in the complex plane of the zeroes of the chromatic polynomial of $C_n$. Unfortunately, we cannot see a direct proof. We can, however, use our recurrence for $\sigma(C_n,x)$ to give a very direct proof of a strong consequence of real zeroes, namely ultra log-concavity of the sequence of coefficients.
\begin{thm} \label{thm-cycles-log-concavity}
For each $n \geq 3$ and $k$ satisfying $1 \leq k \leq n-1$, we have
$$
\left(\frac{S(C_n,k)}{{n \choose k}}\right)^2 \geq \left(\frac{S(C_n,k-1)}{{n \choose k-1}}\right)\left(\frac{S(C_n,k+1)}{{n \choose k+1}}\right).
$$
Moreover for $k \geq \chi(C_n)$ the inequalities above are strict.
\end{thm}
We prove this (together with (\ref{cycle-exp}) and (\ref{cycle-var})) in Section \ref{sec-cycles}, while the remainder of the proof of Theorem \ref{thm-cycles-asy-nor} is in Section \ref{sec-bell}.

\section{Forests} \label{sec-forest}

Let $F^c_n$ be a forest on $n$ vertices with $c$ components (throughout this section we will assume that $c \geq 1$ and that $n \geq c$). It is known that $\sigma(F^c_n,x)$ has all real zeroes; this follows from a general result of Brenti \cite[Theorem 3.20]{Brenti2}, together with Harper's result \cite{Harper} that $\sigma(E_n,x)$ has all real roots. In this section we derive an explicit expression for $S(F^c_n,k)$, which allows us to obtain a number of recurrence relations for $\sigma(F^c_n,x)$. These recurrences allow us to give a direct proof that $\sigma(F^c_n,x)$ has all real zeroes, and moreover allows us to observe patterns between the zeroes of the Stirling polynomials of different forests (Theorem \ref{thm-interlacing}). Our explicit expression for $S(F^c_n,k)$ also allows us to compute the mean and variance of $X^c_n$ (the random variable defined in Theorem \ref{thm-forest-asy-nor}), a key step in the proof of Theorem \ref{thm-forest-asy-nor}.

The chromatic polynomial of $F^c_n$ depends only on $n$ and $c$, and not on $F^c_n$ itself (specifically, if the components of $F_n^c$ have $a_1, \ldots, a_c$ vertices then
$$
\chi_{F^c_n}(x)=\prod_{i=1}^c \left(x(x-1)^{a_i-1}\right) = x^c (x-1)^{n-c}),
$$
and so $S(F^c_n,k)$ also depends only on $n$ and $c$. This allows us to chose a convenient forest to facilitate the calculation of $S(F^c_n,k)$. Let $F^c_n$ consist of a star on $n-c-1$ vertices together with $c-1$ isolated vertices. We have
\begin{equation} \label{clm-forest}
S(F^c_n,k) = \sum_{i \geq 0} {c-1 \choose i} S(n-1-i,k-1).
\end{equation}
Indeed, in a partition counted by $S(F^c_n,k)$, the vertex at the center of the star cannot be in the same class as any of the leaves of the star, but there are no other restrictions. We get a valid partition of the vertex set into $k$ classes by first choosing an arbitrary subset of size $i$ ($i \in \{0, \ldots,c-1\}$) of the isolated vertices to be in the same class as the center of the star, and then choosing an arbitrary partition of the remaining $n-1-i$ vertices into $k-1$ classes. Note that for $c=1$ (the case of trees), the right-hand side of (\ref{clm-forest}) reduces to $S(n-1,k-1)$, a fact observed in \cite{DuncanPeele}.

By (\ref{clm-forest}) we have
\begin{eqnarray}
\nonumber \sigma(F^c_n,x) & = & \sum_{k \geq 0} x^k \sum_{i \geq 0} {c-1 \choose i} S(n-1-i,k-1) \\
\nonumber  & = & x \sum_{i \geq 0} {c-1 \choose i} S_{n-1-i}(x)
\end{eqnarray}
where $S_n(x) = \sum_{k} S(n,k)x^k$ is the generating function of the ordinary Stirling numbers of the second kind. This relation was first observed by Korfhage \cite{Korfhage}, where an inductive proof was given.

It is well known that the polynomials $S_n(x)$ satisfy the recurrence
\begin{equation} \label{rec-stir}
S_n(x) = (x+xD)S_{n-1}(x)
\end{equation}
for $n \geq 2$ (where $D$ is differentiation with respect to $x$), with the initial condition $S_1(x)=x$. It follows by linearity that the polynomials $\sigma(F^c_n,x)$ satisfy
\begin{equation} \label{rec-forest}
\frac{\sigma(F^c_{n+1},x)}{x} = (x+xD)\left(\frac{\sigma(F^c_n,x)}{x}\right)
\end{equation}
for all $n \geq c$. The initial condition for this recurrence is obtained by considering any forest $F_c^c$ on $c$ vertices with $c$ components. There is only one such, the empty graph $E_c$ on $c$ vertices, and so $\sigma(F_c^c,x)=S_c(x)$.

We are now in a position to establish that $\sigma(F^c_n,x)$ has all real zeroes, and moreover that $\sigma(F_{c+1}^c,x) \prec \sigma(F_c^c,x)$ and $\sigma(F_n^c,x) \prec \sigma(F_{n+1}^c,x)$ for all $n \geq c+1$ (items 1 and 2 of Theorem \ref{thm-interlacing}). We begin with the first of these statements.

From \cite{Harper} we know that $\sigma(F_c^c,x)/x$ is a polynomial of degree $c-1$ which is positive at $0$ and whose roots $r_1, \ldots, r_{c-1}$ are all real and satisfy $0 > r_1 > \ldots > r_{c-1}$. Since $\chi(F^c_{c+1},x)=2$, we have that $\sigma(F^c_{c+1},x)$ has a zero of multiplicity $2$ at $0$. We now consider the sign of $\sigma(F^c_{c+1},x)$ at $x=r_i$ for each $i$. By (\ref{rec-forest}) we have
$$
\sigma(F^c_{c+1},r_i) = r_i\sigma(F^c_c,r_i) + r_i^2 D\left.\left(\frac{\sigma(F^c_c,x)}{x}\right)\right|_{x=r_i}.
$$
The first term on the right-hand side above is $0$, and the second has the same sign as that of the derivative of $\sigma(F^c_c,x)/x$ at $x=r_i$. If follows that for odd $i$, $\sigma(F^c_{c+1},r_i)$ is positive, while for even $i$ it is negative. We conclude that $\sigma(F^c_{c+1},x)$ has zeroes between $r_{i+1}$ and $r_i$ for each $i \in \{1, \ldots, c-2\}$. If $c$ is odd then $\lim_{x\rightarrow -\infty} \sigma(F^c_{c+1},x) = \infty$ and so, since $\sigma(F^c_{c+1},r_{c-1})$ is negative, there is a zero of $\sigma(F^c_{c+1},x)$ below $r_{c-1}$, and the same conclusion can be reached by similar reasoning if $c$ is even. This accounts for $c+1$ zeroes of $\sigma(F^c_{c+1},x)$, and since it is a polynomial of degree $c+1$ there are no more. It follows that $\sigma(F_{c+1}^c,x) \prec \sigma(F_c^c,x)$.

For item 2 of Theorem \ref{thm-interlacing} we prove two statements by a parallel induction on $n$. First, for each $n \geq c+1$, the polynomial $\sigma(F^c_n,x)$ has all real non-positive zeroes, has a zero of multiplicity $2$ at $0$, and otherwise has distinct zeroes. Second, for $n \geq c+1$, we have $\sigma(F^c_n,x) \prec \sigma(F^c_{n+1},x)$.
The base case ($n=c+1$) for the first statement has been established in the previous paragraph. For each $n \geq c+1$ we simultaneously establish the second statement for $n$ and the first for $n+1$ using an argument very similar to that presented in the previous paragraph, the details of which we leave to the reader.


The recurrence (\ref{rec-forest}) relates the Stirling polynomials of forests with differing numbers of vertices on the same number of components. There is also a recurrence relating the Stirling polynomials forests on the same number of vertices with different numbers of components, namely
\begin{equation} \label{pascal-forest}
\sigma(F^{c+1}_{n+1},x) = \sigma(F^c_{n+1},x) + \sigma(F^c_n,x)
\end{equation}
To see this, note that
$$
S(F^c_n,k) = \sum_{i \geq 0} {c-1 \choose i} S(n-1-i,k-1) = \sum_{i \geq 0} {c-1 \choose i-1} S(n-i,k-1).
$$
and so, using Pascal's identity,
\begin{eqnarray*}
S(F^c_{n+1},k) + S(F^c_n,k) & = & \sum_{i \geq 0} \left({c-1 \choose i-1} + {c-1 \choose i}\right) S(n-i,k-1) \\
& = & \sum_{i \geq 0} {c \choose i} S(n-i,k-1) \\
& = & S(F^{c+1}_{n+1},k).
\end{eqnarray*}

Using (\ref{pascal-forest}) together with items 1 and 2 of Theorem \ref{thm-interlacing}, we can easily establish items 3, 4 and 5 of the theorem. We begin by considering (\ref{pascal-forest}) when $n=c$. If $r_1 > r_2 \ldots > r_{c-1}$ are the negative zeroes of $\sigma(F^c_c,x)$ and $s_1 > s_2 > \ldots > s_{c-1}$ the negative zeroes of $\sigma(F^c_{c+1},x)$, we know that
$$
0 > r_1 > s_1 > \ldots > r_{c-1} > s_{c-1},
$$
and from (\ref{pascal-forest}) we get that $\sigma(F^{c+1}_{c+1},x)$ is positive at $x=r_i$ and $x=s_i$ for all odd $i$, and negative at $x=r_i$ and $x=s_i$ for all even $i$. Since $\sigma(F^{c+1}_{c+1},x)$ has a zero of multiplicity $1$ at $0$ and has positive derivative at $0$, it follows that $\sigma(F^{c+1}_{c+1},x)$ must have a zero between $r_1$ and $0$, as well as between $r_{i+1}$ and $s_i$ for $i \in \{1, \ldots, c-2\}$. By considering $\lim_{x \rightarrow -\infty} \sigma(F^{c+1}_{c+1},x)$ we also find that there is a zero below $s_{c-1}$. This accounts for all $c+1$ zeroes, and we conclude that $\sigma(F^c_{c+1},x) \prec \sigma(F^{c+1}_{c+1},x)$ and $\sigma(F^c_c,x) \prec \sigma(F^{c+1}_{c+1},x)$. The verification of all the remaining statements in Theorem \ref{thm-interlacing} is very similar, and we leave the details to the reader.

\medskip

We now derive expressions for the mean and variance of $X^c_n$, the number of classes in a uniformly chosen partition of $F^c_n$ into non-empty independent sets, in terms of Bell numbers.
We begin by noting that
\begin{equation} \label{eq-exp_forests}
E(X^c_n) = \frac{\sum_{k \geq 0} k S(F^c_n, k)}{\sum_{k \geq 0} S(F^c_n, k)}.
\end{equation}
Using (\ref{clm-forest}) we see that the denominator of (\ref{eq-exp_forests}) is $\sum_{i \geq 0} {c-1 \choose i} B_{n-1-i}$ (this is the {\em Bell number} of $F^c_n$, the total number of partitions into non-empty independent sets). Using the recurrence $S(n+1,k)=kS(n,k)+S(n,k-1)$ (valid for $(n,k)\neq(0,0)$, with initial conditions $S(0,0)=1$ and $S(n,0)=S(0,k)=0$ for $n, k > 0$), we get
$$
kS(n-1-i,k-1) = S(n-i,k-1) - S(n-1-i,k-2) + S(n-1-i,k-1),
$$
and so, again using (\ref{clm-forest}), the numerator of (\ref{eq-exp_forests}) is $\sum_{i \geq 0} {c-1 \choose i} B_{n-i}$.
This gives (\ref{forest-exp}).

The basic recurrence for the Stirling numbers also gives
\begin{eqnarray*}
k^2S(n-1-i,k-1) & = &  S(n+1-i,k-1) - 2S(n-i,k-2) + 2S(n-i,k-1) \\
&  &  + S(n-1-i,k-3) - 3S(n-1-i,k-2) + S(n-1-i,k-1).
\end{eqnarray*}
From this it follows that
$$
\sum_{i \geq 0} {c-1 \choose i} \sum_{k \geq 0} k^2S(n-1-i,k-1) = \sum_{i \geq 0} {c-1 \choose i} \left(B_{n+1-i}-B_{n-1-i}\right)
$$
and so
$$
E((X^c_n)^2) = \frac{\sum_{i \geq 0} {c-1 \choose i} \left(B_{n+1-i}-B_{n-1-i}\right)}{\sum_{i \geq 0} {c-1 \choose i} B_{n-1-i}},
$$
from which
(\ref{forest-var}) easily follows.

\section{Cycles} \label{sec-cycles}


We begin this section by deriving a recurrence for $S(C_n,k)$ that will allow us to compute the mean and variance of the number of colors in a randomly chosen coloring of the cycle, and which will also be the key to Theorem \ref{thm-cycles-log-concavity}.

Let $P_n$ be the path on $n$ vertices. The graphical Stirling numbers $S(C_n,k)$ satisfy the recurrence
\begin{equation} \label{rec-cycles}
S(C_n,k) = S(n-1,k-1) - S(C_{n-1},k) = S(P_n,k) - S(C_{n-1},k)
\end{equation}
valid for all $n \geq 4$ and all $k$, with initial conditions $S(C_3,3)=1$ and $S(C_3,k)=0$ for $k \neq 3$. To see this recurrence, note that as a special case of (\ref{rec-forest}) we have $S(n-1,k-1)=S(P_n,k)$. Among the partitions of $P_n$ into $k$ non-empty independent sets, those in which the first and last vertices of $P_n$ fall into different classes are in one-to-one correspondence with partitions of the vertex set of $C_n$ into $k$ non-empty independent sets (so there are $S(C_n,k)$ such). The remaining partitions (in which the first and last vertices of $P_n$ fall into the same class) are in one-to-one correspondence with partitions of the vertex set of $C_{n-1}$ into $k$ non-empty independent sets (so there are $S(C_{n-1},k)$ such).

From (\ref{rec-cycles}) we get the polynomial recurrences
\begin{eqnarray}
\label{rec-cycle} \sigma(C_n,x) & = & x S_{n-1}(x) - \sigma(C_{n-1},x) \\
\nonumber & = & \sigma(P_n,x) - \sigma(C_{n-1},x)
\end{eqnarray}
for $n \geq 4$, with initial condition $\sigma(C_3,x)=x^3$. Iterating (\ref{rec-cycle}) we find that
\begin{eqnarray*}
\sigma(C_n,x) & = & x \sum_{i=1}^{n-1} (-1)^{i+1}S_{n-i}(x)\\
&= & \sum_{k \geq 0} x^k \sum_{i \geq 0} (-1)^i S(n-1-i, k-1)
\end{eqnarray*}
a formula that is valid for all $n \geq 2$ if we interpret $C_2$ to be a single edge. Recalling (\ref{rec-stir}), we also have the recurrence
\begin{equation} \label{rec-cycpoly}
\sigma(C_n,x) = x\left((x+xD)\left(\frac{\sigma(C_{n-1},x)}{x}\right)\right) + (-1)^n x^2
\end{equation}
for $n \geq 3$, with initial condition $\sigma(C_2,x)=x^2$.

From here we see no way of obtaining a direct proof that $\sigma(C_n,x)$ has all real zeroes, the problem being that although the operator $x(x+xD)(1/x)$, when applied to a polynomial,  preserves the property of having real zeroes, we cannot say the same with the addition of the $(-1)^n x^2$ term. However, we do easily obtain a slightly weaker statement, the  log-concavity of the sequence $(S(C_n,k)/{n \choose k})_{k=0}^n$ (Theorem \ref{thm-cycles-log-concavity}). Our strategy is to show that the $x(x+xD)(1/x)$ operator preserves ultra log-concavity, and then to deal with the $(-1)^n x^2$ term by observing that it only impacts a bounded number of the ultra log-concavity relations, which can be dealt with by hand using our knowledge of the chromatic polynomial of the cycle.

Let $P(x)=\sum_{k=2}^n a_k x^k$ satisfy $n \geq 3$ and $a_k > 0$ for each $3 \leq k \leq n$. We begin by observing that if $P(x)$ is ultra log-concave then so too is $Q(x)=x(x+xD)(1/x)P(x)$, and that moreover for each $k \geq 3$ (or $k \geq 2$ if $a_2 > 0$) the ultra log-concavity inequalities for $Q(x)$ are all strict. To see this, note that we have
$$
Q(x) = a_2x^2 + \left(\sum_{k=3}^n (a_{k-1} +(k-1)a_k)x^k\right) + a_n x^{n+1}.
$$
To show ultra log-concavity of $Q(x)$ with strict inequalities for $k \geq 3$ (or $k \geq 2$ in the case $a_2>0$), we need to establish the three relations
\begin{equation} \label{inq-ulc1}
\left(\frac{a_2 + 2a_3}{{n+1 \choose 3}}\right)^2 > \left(\frac{a_2}{{n+1 \choose 2}}\right)\left(\frac{a_3 + 3a_4}{{n+1 \choose 4}}\right),
\end{equation}
\begin{equation} \label{inq-ulc2}
\left(\frac{a_{n-1} + (n-1)a_n}{{n+1 \choose n}}\right)^2 > \left(\frac{a_n}{{n+1 \choose n+1}}\right)\left(\frac{a_{n-2} + (n-2)a_{n-1}}{{n+1 \choose n-1}}\right),
\end{equation}
and, for $k=4, \ldots, n-1$,
\begin{equation} \label{inq-ulc3}
\left(\frac{a_{k-1} + (k-1)a_k}{{n+1 \choose k}}\right)^2 > \left(\frac{a_{k-2}+(k-2)a_{k-1}}{{n+1 \choose k-1}}\right)\left(\frac{a_k + ka_{k+1}}{{n+1 \choose k+1}}\right).
\end{equation}
After some algebra, (\ref{inq-ulc1}) is seen to be implied by (thought not equivalent to) $a_3^2/(n-1) \geq a_2a_4/(n-2)$,
which for $n \geq 3$ easily follows from
$\left(a_3/{n \choose 3}\right)^2 \geq \left(a_2/{n \choose 2}\right)\left(a_4/{n \choose 4}\right)$, a consequence of the ultra log-concavity of $P(x)$.
Similarly (\ref{inq-ulc2}) is implied by
$a_{n-1}^2/(n+1) \geq 2a_{n-2}a_n/n$,
a consequence of $\left(a_{n-1}/{n \choose n-1}\right)^2 \geq \left(a_{n-2}/{n \choose n-2}\right)\left(a_n/{n \choose n}\right)$.
This leaves (\ref{inq-ulc3}), which certainly holds if each of
\begin{equation} \label{inq-ulc4}
\left(\frac{a_{k-1}}{{n+1 \choose k}}\right)^2 > \left(\frac{a_{k-2}}{{n+1 \choose k-1}}\right)\left(\frac{a_k}{{n+1 \choose k+1}}\right),
\end{equation}
\begin{equation} \label{inq-ulc5}
(k-1)^2\left(\frac{a_k}{{n+1 \choose k}}\right)^2 > k(k-2)\left(\frac{a_{k-1}}{{n+1 \choose k-1}}\right)\left(\frac{a_{k+1}}{{n+1 \choose k+1}}\right),
\end{equation}
and
\begin{equation} \label{inq-ulc6}
2(k-1)\left(\frac{a_{k-1}}{{n+1 \choose k}}\right)\left(\frac{a_k}{{n+1 \choose k}}\right) -  (k-2)\left(\frac{a_{k-1}}{{n+1 \choose k-1}}\right)\left(\frac{a_k}{{n+1 \choose k+1}}\right) > k\left(\frac{a_{k-2}}{{n+1 \choose k-1}}\right)\left(\frac{a_{k+1}}{{n+1 \choose k+1}}\right)
\end{equation}
hold. From $\left(a_{k-1}/{n \choose k-1}\right)^2 \geq \left(a_{k-2}/{n \choose k-2}\right)\left(a_k/{n \choose k}\right)$ we get (\ref{inq-ulc4})
(after some algebra it reduces to $k^2 > k^2-1$), and from $\left(a_k/{n \choose k}\right)^2 \geq \left(a_{k-1}/{n \choose k-1}\right)\left(a_{k+1}/{n \choose k+1}\right)$ we get (\ref{inq-ulc5}) (it reduces to $(k-1)^2(n-k+1)^2 > ((k-1)^2-1)((n-k+1)^2-1)$). Finally (\ref{inq-ulc6}) follows (after some algebra) for all $k$ in the given range (and more generally for $n > k-2$) from another consequence of the ultra log-concavity of $P(x)$, namely $\left(a_{k-1}/{n \choose k-1}\right)\left(a_k/{n \choose k}\right) \geq \left(a_{k-2}/{n \choose k-2}\right)\left(a_{k+1}/{n \choose k+1}\right)$.

We now prove Theorem \ref{thm-cycles-log-concavity} by induction on $n$. The base case $n = 3$ is trivial. If $n \geq 3$ is odd, then (by the inductive hypothesis and the observation of the previous paragraphs), $x(x+xD)(\sigma(C_{n-1}, x)/x)$ is ultra log-concave with strict inequalities for $k \geq 2$, and, since $C_{n-1}$ is connected and bipartite, the coefficient of $x^2$ is $1$. It easily follows that $\sigma(C_n,x)=x(x+xD)(\sigma(C_{n-1}, x)/x)-x^2$ is ultra log-concave with strict inequalities for $k \geq 3$. For even $n\geq 3$, $x(x+xD)(\sigma(C_{n-1}, x)/x)$ is also ultra log-concave with strict inequalities for $k \geq 3$, but, since $C_{n-1}$ is not bipartite, the coefficient of $x^2$ is $0$. To conclude ultra log-concavity of $\sigma(C_n, x) = x(x+xD)(\sigma(C_{n-1}, x)/x) + x^2$ with strict inequalities for $k \geq 3$ we also need to show that
\begin{equation} \label{log-con-to-show}
\left(\frac{S(C_n, 3)}{{n \choose 3}}\right)^2 > \frac{S(C_n,4)}{{n \choose 4}{n \choose 2}}.
\end{equation}
Using $\chi_{C_n}(x) = (-1)^n(x-1) + (x-1)^n$ and (\ref{eq-chrom}) we get
$$
S(C_n, 3) = \frac{\chi_{C_n}(3) - 3\chi_{C_n}(2)}{3!} = \frac{2^n - (-1)^n - 3}{6}
$$
and
$$
S(C_n, 4) = \frac{\chi_{C_n}(4) - 4\chi_{C_n}(3) + 6\chi_{C_n}(2)}{4!} = \frac{3^n - 4\cdot2^n + (-1)^n + 6}{24}
$$
from which (\ref{log-con-to-show}) follows for $n \geq 4$.

\medskip

Our recurrence for $S(C_n,k)$ also allows us to compute the expectation and variance of $X^{\rm cycle}_n$, the number of classes in a uniformly chosen partition of the vertex set of $C_n$ into non-empty independent sets. Iterating (\ref{rec-cycles}) we get
$$
S(C_n, k) = \sum_{i \geq 0} (-1)^i S(n-1-i,k-1)
$$
and so
$$
E(X^{\rm cycle}_n) = \frac{\sum_{i \geq 0}  (-1)^i \sum_{k \geq 0} k S(n-1-i,k-1)}{\sum_{i \geq 0} (-1)^i \sum_{k \geq 0} S(n-1-i,k-1)}
$$
and
$$
E((X^{\rm cycle}_n)^2) = \frac{\sum_{i \geq 0}  (-1)^i \sum_{k \geq 0} k^2 S(n-1-i,k-1)}{\sum_{i \geq 0} (-1)^i \sum_{k \geq 0} S(n-1-i,k-1)}.
$$
The  calculations made in Section \ref{sec-forest} to verify (\ref{forest-exp}) and (\ref{forest-var}) also easily yield (\ref{cycle-exp}) and (\ref{cycle-var}).

\section{Asymptotic normality} \label{sec-asy-nor}

The connection between real zeros and asymptotic normality was first made by L\'evy \cite{Levy}, and rediscovered by Harper \cite{Harper}, who was considering the specific problem of the Stirling numbers of the second kind. It has since been presented in a more general form by numerous authors, including Bender \cite{Bender}, Canfield \cite{Canfield} and Godsil \cite{Godsil2}. The following general result is from \cite{Bender}.
\begin{thm} \label{thm-asy-nr}
Let $({\mathcal C}_n)_{n \geq 0}$ be a sequence of sets, and let each ${\mathcal C}_n$ come with a partition ${\mathcal C}_n = \cup_{k=0}^n {\mathcal C}_n^k$. Set $|{\mathcal C}_n|=c_n$ and $|{\mathcal C}^k_n|=c^k_n$. Suppose that for each $n$ the polynomial $P_n(x) = \sum_{k=0}^n c_n^k x^k$ has all real zeroes. Let $X_n$ be the upper index of a uniformly chosen element $\omega$ of ${\mathcal C}_n$ (that is, $X_n(\omega)=k$ if and only if $\omega \in {\mathcal C}_n^k$). If ${\rm Var}(X_n) \rightarrow \infty$ as $n \rightarrow \infty$ then $(X_n)_{n \geq 0}$ is asymptotically normal. Moreover, we have the {\em local limit theorem}
$$
\lim_{n \rightarrow \infty} \sup_{x \in {\mathbb R}} \left|\sqrt{{\rm Var}(X_n)}\Pr\left(X_n=\left[E(X_n)+x\sqrt{{\rm Var}(X_n)}\right]\right)  - \frac{1}{\sqrt{2\pi}}e^{-x^2/2}  \right| = 0.
$$
\end{thm}

The advantage of a local limit theorem is that it provides quantitative information about the $c_n^k$'s that asymptotical normality does not; see for example \cite{Bender}, \cite{Canfield} or \cite{CzabarkaErdosJohnsonKupczokSzekely} for details.

For completeness we provide a short proof of the asymptotic normality part of Theorem \ref{thm-asy-nr} that gives explicit information about the rate of convergence. The probability generating function of $X_n$ is $\Phi_n(x)=P_n(x)/P_n(1)$. Since $\Phi_n(x)$ has all real zeroes and all non-negative coefficients, we may factor it as
$$
\Phi_n(x)=\prod_{i=1}^n \left(\frac{\lambda_i}{1+\lambda_i} + \frac{x}{1+\lambda_i}\right)
$$
with each $\lambda_i$ non-negative. It follows that we may write
$$
X_n=\sum_{i=0}^n X_n^{(i)}
$$
where the $X_n^{(i)}$'s are independent, and each $X_n^{(i)}$ takes value $0$ with probability $\lambda_i/(1+\lambda_i)$ and value $1$ with probability $1/(1+\lambda_i)$, and so
$$
X_n - E(X_n) = \sum_{i=1}^n Y_i
$$
where the $Y_i$'s are a collection of independent, mean $0$ random variables, with specifically $Y_i$ taking value $-1/(1+\lambda_i)$ with probability $\lambda_i/(1+\lambda_i)$ and value $\lambda_i/(1+\lambda_i)$ with probability $1/(1+\lambda_i)$. Since $Y_i$ only takes values between $-1$ and $1$, we have $E(|Y_i|^3) \leq E(Y_i^2)$.

The Berry-Esseen theorem \cite{Berry} now says that there is an absolute constant $C>0$ such that for all $x \in {\mathbb R}$ we have
$$
\left|\Pr\left(\frac{\sum_{i=1}^n Y_i}{\sqrt{{\rm Var}\left(\sum_{i=1}^n Y_i\right)}} \leq x\right) - \frac{1}{\sqrt{2\pi}} \int_{-\infty}^x e^{-x^2/2}~\!dx\right| \leq \frac{C}{\sqrt{{\rm Var}(X_n)}},
$$
from which asymptotic normality follows.

\section{Proofs of Theorems \ref{thm-forest-asy-nor} and \ref{thm-cycles-asy-nor}} \label{sec-bell}

We begin by establishing (\ref{forest-exp-est}) and (\ref{forest-var-est}). We need two preliminary estimates involving Bell numbers. The first is due to Canfield and Harper \cite{CanfieldHarper}. We have
\begin{equation} \label{inq-knuth}
\frac{B_{n-1}}{B_n} = \frac{W(n)}{n}\left(1+O\left(\frac{1}{n}\right)\right) = \frac{W(n)}{n} + O\left(\frac{\log n}{n^2}\right)
\end{equation}
(where recall $W(n) = \log n -\log \log n + O(1)$ is the Lambert $W$ function, the unique positive solution to $W(n)e^{W(n)} = n$). The second is due to Harper \cite{Harper}. We have
\begin{equation} \label{inq-harper}
\frac{B_{n+2}}{B_n} - \left(\frac{B_{n+1}}{B_n}\right)^2 = \frac{n}{W(n)(W(n)+1)} + O\left(\frac{1}{\log^2 n}\right).
\end{equation}
(Harper gives an error term of $o(1)$; see for example \cite{CzabarkaErdosJohnsonKupczokSzekely} for the more refined error.)

We need a slight strengthening of (\ref{inq-knuth}). For $k=o(n)$, (\ref{inq-knuth}) gives
\begin{equation} \label{eq-usingKnuth}
\frac{B_{n-k-1}}{B_{n-k}} = \frac{W(n-k)}{n-k} + O\left(\frac{\log (n-k)}{(n-k)^2}\right) = \frac{W(n)}{n} + O\left(\frac{k\log n}{n^2}\right),
\end{equation}
with the same bound for $B_{n+k-1}/B_{n+k}$. Here we have used the mean value theorem to estimate $W(n)-W(n-k)$, noting that $W'(x)=W(x)/(x(1+W(x)))$. A similar argument gives
\begin{equation} \label{inq-harper2}
\frac{B_{n+1}}{B_{n-1}} - \left(\frac{B_n}{B_{n-1}}\right)^2 = \frac{n}{W(n)(W(n)+1)} + O\left(\frac{1}{\log^2 n}\right).
\end{equation}

From (\ref{eq-usingKnuth}) we get that for $m \in \{n-1, n, n+1\}$ we have
\begin{eqnarray*}
\sum_{i \geq 0} {c-1 \choose i}B_{m-i} & = & B_m\left(1+\frac{W(n)}{n} + O\left(\frac{c\log n}{n^2}\right)\right)^c \\
& = & B_m\left(1+\frac{W(n)}{n}\right)\left(1+O\left(\frac{c^2\log n}{n^2}\right)\right)
\end{eqnarray*}
(the second equality valid for $c=o(n/\sqrt{\log n})$).

It now immediately follows that $E(X_n^c) = (n/W(n))(1+O(1/n))$ (recall that we are in the range $c=O(\sqrt{n/\log n}))$. It also follows that
\begin{eqnarray*}
\frac{\sum_{i=0}^{c-1} {c-1 \choose i}B_{n+1-i}}{\sum_{i=0}^{c-1} {c-1 \choose i}B_{n-1-i}} - \left(\frac{\sum_{i=0}^{c-1} {c-1 \choose i}B_{n-i}}{\sum_{i=0}^{c-1} {c-1 \choose i}B_{n-1-i}}\right)^2 & = & \frac{B_{n+1}}{B_{n-1}} - \left(\frac{B_n}{B_{n-1}}\right)^2 + O\left(\frac{c^2}{\log n}\right) \\
& = & \frac{n}{W(n)(W(n)+1)} + O\left(\frac{c^2}{\log n}\right)
\end{eqnarray*}
(here using (\ref{eq-usingKnuth}) and (\ref{inq-harper2})). As long as $c < C\sqrt{n/\log n}$ for some sufficiently small $C>0$, we have ${\rm Var}(X_n)=\omega(1)$, and so we complete the proof of Theorem \ref{thm-forest-asy-nor} by appealing to Theorem \ref{thm-asy-nr}.

\medskip

We now turn to Theorem \ref{thm-cycles-asy-nor}. For each $k \in \{-1, 0, 1\}$ we use (\ref{eq-usingKnuth}) to obtain
\begin{eqnarray*}
\sum_{i \geq 0} (-1)^i B_{n+k-i} & = & B_{n+k} \left(1-\frac{B_{n+k-1}}{B_{n+k}} + \frac{B_{n+k-2}B_{n+k-1}}{B_{n+k-1}B_{n+k}} + \right. \\
& & ~~~~~\left.\frac{B_{n+k-3}B_{n+k-2}B_{n+k-1}}{B_{n+k-2}B_{n+k-1}B_{n+k}} + O\left(\frac{n B_{n+k-4}}{B_{n+k}}\right)\right) \\
& = & B_{n+k}\left(1-\frac{W(n)}{n} + \frac{W(n)^2}{n^2} -  \frac{W(n)^3}{n^3}\right)\left(1+O\left(\frac{\log n}{n^2}\right)\right).
\end{eqnarray*}
Using (\ref{cycle-exp}) and (\ref{cycle-var}) this gives
$$
E(X_n^{\rm cycle}) = \frac{B_n}{B_{n-1}}\left(1+O\left(\frac{\log n}{n^2}\right)\right) = \frac{n}{W(n)} + O\left(\frac{1}{\log n}\right)
$$
and
\begin{eqnarray*}
{\rm Var}(X^{\rm cycle}_n) & = & \frac{B_{n+1}}{B_{n-1}} - \left(\frac{B_n}{B_{n-1}}\right)^2 + O\left(\frac{1}{\log n}\right) \\
& = & \frac{n}{W(n)(W(n)+1)} + O\left(\frac{1}{\log n}\right),
\end{eqnarray*}
the last equality from (\ref{inq-harper2}). Since ${\rm Var}(X^{\rm cycle}_n) = \omega(1)$, we complete the proof of Theorem \ref{thm-cycles-asy-nor} by appealing to Theorem \ref{thm-asy-nr}.

\end{document}